\documentclass{article}
\usepackage{amsmath}
\usepackage{amsfonts}
\usepackage{amsthm}
\usepackage[all]{xy}
\newcommand{\ZZ}{\mathbb{Z}}
\newtheorem{theorem}{Theorem}
\newtheorem{lemma}{Lemma}
\newtheorem{cor}{Corollary}
\newtheorem{prop}{Proposition}
\theoremstyle{definition}
\newtheorem{defn}{Definition}
\newtheorem{example}{Example}
\theoremstyle{remark}
\newtheorem*{acn}{Acnowledgements}
\DeclareMathOperator{\Hom}{Hom}
\DeclareMathOperator{\har}{char}
\DeclareMathOperator{\inv}{inv}
\DeclareMathOperator{\symm}{sym}
\DeclareMathOperator{\spa}{span}
\begin{document}
\title{Diagrams of representations}

\author{Aleksandrs Mihailovs\\
Department of Mathematics\\
University of Pennsylvania\\
Philadelphia, PA 19104-6395\\
mihailov@math.upenn.edu\\
http://www.math.upenn.edu/$\sim$mihailov/
}
\date{\today}
\maketitle

\begin{abstract}
For a representation of a Lie algebra, one can construct a diagram 
of the representation, i.\ e.\ a directed graph with edges labeled by matrix 
elements of the representation. This article explains how to use 
these diagrams to describe normal forms, orbits and invariants of 
the representation, especially for the case of nilpotent Lie algebras.
\end{abstract}    
\setlength{\baselineskip}{1.5\baselineskip}

\section{Introduction}\label{sec0}

For a representation of a Lie algebra, one can construct a diagram 
of the representation, i.\ e.\ a directed graph with edges labeled by matrix 
elements of the representation. In section \ref{sec1} I give the detailed 
definition and a few examples. In section \ref{sec2} I describe how to construct 
the diagrams of symmetric and exterior powers of representations, with a few 
examples. Section \ref{sec3} is the main section of this work. I introduce there  
the diagram method for the description of normal forms on the orbits of a strictly 
triangular representation of a Lie algebra, usually nilpotent or pronilpotent. 
Some results in the beginning of this section (without using the diagrams) were 
obtained independently and earlier by Igor Brodski \cite{Br}. Then, in 
section \ref{sec4}, I show how to apply the diagram method to the description of 
normal forms of quadratic differentials on a line. First these normal forms were 
found by Alexandre Kirillov \cite{KV}. In section \ref{sec5} I give the analogous 
description of normal forms of generalized tensor fields on a line, using the 
diagram method. The last two sections are closely connected to my previous work 
\cite{MF}.

\section{The diagrams of representations of Lie algebras}\label{sec1} 

Let $f$ be a field of characteristic $p\geq 0$, $L$ a Lie $f$-algebra and 
$E$ an $f$ space. 
Denote $L^\ast=\Hom_f(L,f)$ and $E^\ast=\Hom_f(E,f)$, the $f$-spaces of 
$f$-linear forms on $L$ and $E$. Fix a basis $(e_i)_{i\in I}$ of $E$ and the 
dual basis $(x_i)_{i\in I}$ of $E^\ast$.  

For a representation $T$ of a Lie algebra $L$ in 
$E$, define an $f$-bilinear mapping 
\begin{equation}
\phi_T: E^\ast\times E\rightarrow L^\ast,\quad \phi_T(x,e)l=x(T(l)e) .
\end{equation}

\begin{defn}\label{defn1}
A diagram of a representation $T$ of a Lie algebra $L$ in $E$, 
corresponding to a basis $(e_i)_{i\in I}$ of $E$ is 
a directed graph with the set of vertices $I$, having an edge $(i,j)$ iff 
$\phi_T(x_j, e_i)\neq 0$, in which case this edge is labeled by 
$\phi_T(x_j, e_i)\in L^\ast$. We suppose here that this graph doesn't have multiple 
edges, but can have loops. $(x_j)_{j\in I}$ denotes a basis of $E^\ast$, dual to 
$(e_i)_{i\in I}$, as usual.
\end{defn}

\begin{example}\label{ex1}
For the adjoint representation of a Heisenberg Lie algebra $H$ with generators 
$X, Y, Z$ satisfying 
\begin{equation}
[X,Y]=Z,\quad [X,Z]=[Y,Z]=0 ,
\end{equation}
the diagram is \\ 
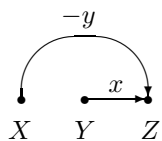
\begin{picture}(410, 55) 
\put(187,18){\circle*{3}}
\put(211,18){\circle*{3}}
\put(235,18){\circle*{3}}
\put(182,3){$X$}
\put(207,3){$Y$}
\put(232,3){$Z$}
\put(212,18){\vector(1,0){22}}
\put(211,18){\oval(48,48)[t]}
\put(235,19){\vector(0,-1){0}}
\put(220,21){$x$}
\put(202,47){$-y$}
\end{picture}\\
where $(x,y,z)$ is the basis of $H^\ast$, dual to $(X,Y,Z)$.
\end{example}

The diagrams are very useful for the description of various operations on 
representations. For instance, 

\begin{prop}\label{prop1}
The diagram of the dual representation $T^\ast$  
corresponding to the basis $(x_j)_{j\in I}$ of $E^\ast$, can be obtained 
from the diagram of the representation $T$ corresponding to the dual basis 
$(e_i)_{i\in I}$ of $E$, by  
changing the directions of all arrows and all the signs of their labels. 
\end{prop}
\begin{proof}By definition.\end{proof}

\begin{example}\label{ex2}
The diagram of the co-adjoint representation of $H$, see example \ref{ex1}, is 

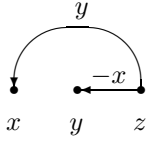
\begin{picture}(410, 55) 
\put(187,18){\circle*{3}}
\put(211,18){\circle*{3}}
\put(235,18){\circle*{3}}
\put(184,3){$x$}
\put(208,3){$y$}
\put(232,3){$z$}
\put(234,18){\vector(-1,0){22}}
\put(211,18){\oval(48,48)[t]}
\put(187,19){\vector(0,-1){0}}
\put(216,21){$-x$}
\put(210,47){$y$}
\end{picture}.
\end{example}

The underlying graph of a diagram contains important information about the 
representation. For the adjoint representation, it contains a lot of information 
about the Lie algebra. In particular, 

\begin{prop}\label{prop2} 
Let $G$ be the underlying directed graph of the adjoint representation of 
a Lie algebra $L$. If $G$ is a graph without edges, $L$ is 
a commutative Lie algebra. If $G$ is a graph without oriented cycles, 
$L$ is a nilpotent Lie algebra of class of nilpotency $\leq l(G)+1$ where 
$l(G)$ is the 
length of $G$, i.\ e.\ the number of edges in the longest directed path 
which is a subgraph of $G$.  
\end{prop}
\begin{proof}Again, everything follows directly from the definitions. \end{proof}

\begin{example}\label{ex3}
For the Heisenberg Lie algebra $H$, see example \ref{ex1}, $l(G)=1$. 
By Proposition \ref{prop2}, $H$ is a nilpotent Lie algebra of class $\leq 2$, i.\ 
e.\ , a metabelian Lie algebra.
\end{example}

\begin{defn}\label{defn2}
If the underlying graph of a diagram $C$ is a subgraph of the underlying graph 
of a diagram $D$ containing with each vertex all the edges starting from this  
vertex, with the same labels, we'll say that the diagram $C$ is a subdiagram of 
the diagram $D$. If the underlying graph of a diagram $C$ is a subgraph of the 
underlying graph of a diagram $D$ containing with each vertex all the edges 
ending in this vertex, with the same labels, we'll say that the diagram $C$ is 
a quotient diagram of the diagram $D$.
\end{defn}

\begin{prop}\label{prop3}
If $D$ is a diagram 
of a representation $T$ of a Lie algebra $L$ in an $f$-space $E$ corresponding to 
the basis $(e_i)_{i\in I}$ and $C$ is a subdiagram of $D$, 
then the $f$-linear subspace $E_C$ 
spanned by $(e_i)_{i\in J}$ where $J$ is the set of the vertices of $C$, is 
$T$-invariant and $C$ is the 
diagram of the restriction of $T$ on $E_C$, 
corresponding to the basis $(e_i)_{i\in J}$. Moreover, if we 
denote $\overline{C}$ the diagram with the set of vertices $I\setminus J$, containing   
all the edges of $D$ between them, with the same labels, then $\overline{C}$ is a 
quotient diagram coinciding with the diagram of the quotient representation of 
$L$ in the quotient space $E/E_C$, corresponding to the basis $(e_i+E_C)_{i\in 
I\setminus J}$. 
\end{prop} 
\begin{proof}By definition.\end{proof}

Dually,
\begin{prop}\label{prop4}
If $D$ is a diagram of a representation $T$ of a Lie algebra $L$ in an $f$-space 
$E$ corresponding to the basis $(e_i)_{i\in I}$ and $C$ is a quotient diagram of 
$D$, then the diagram $\overline{C}$ defined the same as in Proposition \ref{prop3}, 
as a complete subgraph of the underlying graph of $D$ with the set of vertices 
$I\setminus J$ where $J$ is the set of vertices of $C$, with the same labeling of 
the edges, is a subdiagram of $D$ and $\overline{\overline{C}}
=C$.
\end{prop}
\begin{proof}Again, it follows directly from the definitions.\end{proof}
    
\begin{prop}\label{prop5}
Let $T$ and $U$ be representations of a Lie algebra $L$ in $f$-spaces $V$ and $W$ 
with bases $(v_i)_{i\in I}$ and $(w_j)_{j\in J}$ respectively. Denote 
$D_T,D_U,D_{T\oplus U},D_{T\otimes U}$ the diagrams of the representations 
$T,U,T\oplus U, T\otimes U$ corresponding to the bases $(v_i)_{i\in I}$, 
$(w_j)_{j\in J}$, $(v_i)_{i\in I}\amalg (w_j)_{j\in J}$ and $(v_i\otimes w_j)_
{(i,j)\in I\times J}$. Then 
\begin{gather}\label{pr5eq1}
D_{T\oplus U}=D_T\amalg D_U \\
D_{T\otimes U}=D_T\times D_U , \label{pr5eq2}
\end{gather}
meaning that $D_T\amalg D_U$ is a disjoint union of $D_T$ and $D_U$ and 
$D_T\times D_U$ is a reduced Cartesian product of $D_T$ and $D_U$ (`reduced' means 
that we 
replace multiple loops in a vertix of the Cartesian product by one loop labeled 
by the sum of the labels of these multiple loops).
\end{prop}
\begin{proof}It is evident.
\end{proof}

At the end of the section, note that 
the diagram contains all the information about the representation: 

\begin{prop}\label{prop6}
If $D$ is a diagram of a representation $T$ of a Lie algebra $L$ in $V$ 
corresponding to the basis $(e_i)_{i\in I}$, then 
\begin{equation}\label{eqpr6}
T(l)\sum_{i\in I}x_ie_i=\sum_{(i,j)\in E(D)}x_iw(i,j)(l)e_j  
\end{equation}
where $E(D)$ is the set of edges of the underlying graph of the diagram $D$, 
and $w(i,j)\in L^\ast$ denotes the label of the edge $(i,j)$.
\end{prop}
\begin{proof}It follows from Definition \ref{defn1}.\end{proof} 

In other words, 
\begin{cor}\label{cor1}
If $D$ is a diagram of a representation $T$ of a Lie algebra $L$ in $V$ 
corresponding to the basis $(e_i)_{i\in I}$, then for any $l\in L$, the 
matrix of $T(l)$ in the basis $(e_i)_{i\in I}$ is $^tW(D)(l)$, the transposed 
weight matrix of the weighted graph $D$, applied to $l$ supposing that for 
the matrix elements $a_{ij}\in L^\ast$ one has 
\begin{equation}\label{c1eq1}
\begin{pmatrix}a_{11}&\dots &a_{1n}\\
\hdotsfor[1.5]{3}\\ 
a_{n1}&\dots &a_{nn}\end{pmatrix}(l)= 
\begin{pmatrix}a_{11}(l)&\dots &a_{1n}(l)\\
\hdotsfor[2]{3}\\ 
a_{n1}(l)&\dots &a_{nn}(l)\end{pmatrix}
\end{equation}
\end{cor}
\begin{proof}It follows from Proposition \ref{prop6}.\end{proof} 

\section{Diagrams of $\lambda$-operations on representations}\label{sec2}
  
Proposition \ref{prop5} shows that the diagrams are useful for the description 
of sums and tensor products of representations. They are extremely useful for 
the description of other operations as well.

\begin{defn}\label{defn3}
For a reduced weighted directed graph $G$
with linearly ordered set of vertices $V(G)$, denote $S^nG$ a reduced 
weighted directed graph with the set of vertices 
\begin{equation}\label{d3eq1} 
V(S^nG)=S^nV(G)=\{(v_1,\dots,v_n)\in V(G)^n\quad |\quad v_1\leq\dots\leq v_n\} .
\end{equation}
For $i=1,\dots,n$ denote $\Pr_i:V(S^nG)\rightarrow V(S^{n-1}G)$ the projection 
obtained by the erasing of the vertices $v_i$ standing on the $i$-th place of the 
sequence \eqref{d3eq1}.   
For $A,B\in V(S^nG)$ so that $A\neq B$, the edge $(A,B)\in E(S^nG)$ exists iff 
for some $i,j$ one has $Pr_i(A)=Pr_j(B)$, there is an edge 
$(v_i(A),v_j(B))\in E(G)$ where $v_i$ denotes the vertex of $G$ standing on the 
$i$-th place of the sequence \eqref{d3eq1} and 
\begin{equation}\label{d3eq2}
(\text{$\#$ of the entries of $v_i(A)$ in $A$})\cdot w(v_i(A),v_j(B))\neq 0
\end{equation}
where $w$ denotes the weight function; 
in which case the weight of the edge $(A,B)\in E(S^nG)$ equals \eqref{d3eq2}.
For $A\in V(S^nG)$, the edge $(A,A)\in E(S^nG)$ exists iff 
\begin{equation}\label{d3eq3}
\sum_{i=1}^nw(v_i(A),v_i(A)) \neq 0 ,
\end{equation}
in which case the weight of the edge $(A,A)\in E(S^nG)$ equals \eqref{d3eq3}.
\end{defn}

\begin{prop}\label{prop7}
If $D$ is a diagram of a representation $T$ of a Lie algebra $L$ in $V$ 
corresponding to the basis $(e_i)_{i\in I}$, then $S^nD$ is a diagram of 
the representation $S^nT$ of $L$ in $S^nV=V^{\otimes n}/{\cal I}_n$ where 
${\cal I}_n$ is a subspace of $V^{\otimes n}$ spanned by 
$(e_\alpha-e_{\alpha^s})_{\alpha\in I^n, s\in S_n}$ where $e_{i_1\dots i_n} = 
e_{i_1}\otimes\dots\otimes e_{i_n}$, corresponding to the basis 
$(e_{i_1}\dots e_{i_n})_{(i_1,\dots,i_n)\in S^nI}$ where $e_{i_1}\dots e_{i_n}$ 
is the projection of $e_{i_1\dots i_n}$.
\end{prop} 
\begin{proof}It follows directly from the definitions. \end{proof}

\begin{example}\label{ex5}
Let $T=\rho_1$ be the standard $2$-dimensional representation of 
$sl(2)$, the Lie algebra of $2\times 2$ matrices with trace 0. Denote 
\begin{equation}\label{ex5eq1}
X_+=\begin{pmatrix}0&0\\ 1&0\end{pmatrix},\quad 
X_-=\begin{pmatrix}0&-1\\ 0&0\end{pmatrix},\quad 
H=\begin{pmatrix}1&0\\ 0&-1\end{pmatrix}
\end{equation}
the standard basis elements of $sl(2)$ and $(x_+, x_-, h)$ the dual basis of 
$sl(2)^\ast$. The diagram of $T$ corresponding to the standard basis 
$(u=\binom{1}{0}, v=\binom{0}{1})$, is 
\begin{equation}\label{ex5eq2}
\xymatrix{u \ar@(ul,ur)[]^h \ar@/^/[r]^{x_+} 
&v \ar@(ul,ur)[]^{-h} \ar@/^/[l]^{-x_-}} .
\end{equation}
It follows from Proposition \ref{prop7}, that the diagram for $S^2T=\rho_2$ is
\begin{equation}\label{ex5eq3}
\xymatrix{u^2 \ar@(ul,ur)[]^{2h} \ar@/^/[r]^{2x_+} 
&uv \ar@/^/[l]^{-x_-} \ar@/^/[r]^{x_+} 
&v^2 \ar@(ul,ur)[]^{-2h} \ar@/^/[l]^{-2x_-}} ,
\end{equation}
and the diagram for $S^3T=\rho_3$ is
\begin{equation}\label{ex5eq4}
\xymatrix{u^3 \ar@(ul,ur)[]^{3h} \ar@/^/[r]^{3x_+}
&u^2v \ar@/^/[l]^{-x_-} \ar@(ul,ur)[]^h \ar@/^/[r]^{2x_+}
&uv^2 \ar@/^/[l]^{-2x_-} \ar@(ul,ur)[]^{-h} \ar@/^/[r]^{x_+}
&v^3 \ar@(ul,ur)[]^{-3h} \ar@/^/[l]^{-3x_-}} .
\end{equation}
We supposed in \eqref{ex5eq3} and \eqref{ex5eq4} that $p=\har f\neq 2$ and $\neq 3$. 
For $p=2$ or $p=3$ one has to delete from the diagrams \eqref{ex5eq3} and 
\eqref{ex5eq4} all the edges labeled $\pm px_\pm$ and $\pm ph$.
\end{example}

\begin{defn}\label{defn4}
For a reduced weighted directed graph $G$ with linearly ordered set of vertices $V(G)$ denote 
$\Lambda^nG$ a reduced weighted directed graph with the set of vertices  
\begin{equation}\label{d4eq1}
V(\Lambda^nG)=\Lambda^nV(G)=\{(v_1,\dots,v_n)\in V(G)^n\quad |\quad v_1<\dots< v_n\}
\end{equation}
For $i=1,\dots,n$ denote $\Pr_i:V(\Lambda^nG)\rightarrow V(\Lambda^{n-1}G)$ the projection 
obtained by the erasing of the vertices $v_i$ standing on the $i$-th place of the 
sequence \eqref{d4eq1}.
For $A,B\in V(\Lambda^nG)$ so that $A\neq B$, 
the edge $(A,B)\in E(\Lambda^nG)$ exists iff 
for some $i,j$ one has $Pr_i(A)=Pr_j(B)$ and there is an edge 
$(v_i(A),v_j(B))\in E(G)$ where $v_i$ denotes the vertex of $G$ standing on the 
$i$-th place of the sequence \eqref{d4eq1}, in which case the weight 
of the edge $(A,B)\in E(\Lambda^nG)$ is 
\begin{equation}\label{d4eq2}
w(A,B)=(-1)^{j-i}w(v_i,v_j)
\end{equation}  
where $w$ denotes the weight. For 
$A\in V(\Lambda^nG)$, there is an edge $(A,A)\in E(\Lambda^nG)$ iff 
\begin{equation}\label{d4eq3}
\sum_{i\in A}w(i,i) \neq 0 ,
\end{equation}
in which case the weight of the edge $(A,A)\in E(\Lambda^nG)$ equals \eqref{d4eq3}.
\end{defn}

\begin{prop}\label{prop8}
If $D$ is a diagram of a representation $T$ of a Lie algebra $L$ in $V$ 
corresponding to the basis $(e_i)_{i\in I}$, then $\Lambda^nD$ is a diagram of 
the representation $\Lambda^nT$ of $L$ in $\Lambda^nV=V^{\otimes n}/{\cal J}_n$ 
where ${\cal J}_n$ is a subspace of $V^{\otimes n}$ spanned by 
such tensor products $e_{i_1}\otimes\dots\otimes e_{i_n}$ that 
$i_k=i_l$ for some $k\neq l$; corresponding to the basis 
$(e_{i_1}\wedge\dots\wedge e_{i_n})_{(i_1<\dots<i_n)\in \Lambda^nI}$ 
where $e_{i_1}\wedge\dots\wedge e_{i_n}$ 
is the projection of $e_{i_1}\otimes\dots\otimes e_{i_n}$.
\end{prop} 
\begin{proof}Again, it follows directly from the definitions. \end{proof}

Sometimes one has to consider a subrepresentations $S_nT$ and $\Lambda_nT$ instead 
of the corresponding quotient representations described above. The diagrams of them 
can be easily described as well. 

\begin{defn}\label{defn5}
For a reduced weighted directed graph $G$ 
with a linearly ordered set of vertices $V(G)$, denote $S_nG$ a reduced 
weighted directed graph with the same set of vertices $V(S_nG)=V(S^nG)$ as 
$S^nG$, see \eqref{d3eq1}. The same as for $S^nG$, 
for $A\in V(S_nG)$, the edge $(A,A)\in E(S_nG)$ exists iff
\begin{equation}\label{d5eq1}
\sum_{i=1}^nw(v_i(A),v_i(A)) \neq 0 ,
\end{equation}
in which case the weight of the edge $(A,A)\in E(S_nG)$ equals \eqref{d5eq1}. 
The difference between Definition \ref{defn3} and this definition is 
that for $A,B\in V(S_nG)$ so that $A\neq B$, the edge $(A,B)\in E(S_nG)$ exists iff 
for some $i,j$ one has $Pr_i(A)=Pr_j(B)$, there is an edge 
$(v_i(A),v_j(B))\in E(G)$ and 
\begin{equation}\label{d5eq2}
(\text{$\#$ of the entries of $v_j(B)$ in $B$})\cdot w(v_i(A),v_j(B))\neq 0 ,
\end{equation}
in which case the weight of the edge $(A,B)\in E(S_nG)$ equals \eqref{d5eq2}.
\end{defn}

\begin{prop}\label{prop9}
If $D$ is a diagram of a representation $T$ of a Lie algebra $L$ in $V$ 
corresponding to the basis $(e_i)_{i\in I}$, then $S_nD$ is a diagram of 
the subrepresentation $S_nT$ of $L$ in $S_nV\subseteq V^{\otimes n}$ where 
$S_nV$ is a subspace of $V^{\otimes n}$ with the basis $(e_\alpha^{\symm})_
{\alpha\in S^nI}$ where for $\alpha\in I^n$ 
\begin{equation}\label{p9eq1}
e_\alpha^{\symm}=\sum_{\beta\in S_n(\alpha)}e_\beta
\end{equation}
with $e_{i_1\dots i_n} = 
e_{i_1}\otimes\dots\otimes e_{i_n}$; corresponding to this basis. 
\end{prop} 
\begin{proof}It follows directly from the definitions. \end{proof}

\begin{example}\label{ex4}
\begin{gather}
e_{111}^{\symm}=e_1\otimes e_1\otimes e_1 ,\\
e_{112}^{\symm}=e_1\otimes e_1\otimes e_2+e_1\otimes e_2\otimes e_1+
e_2\otimes e_1\otimes e_1 .
\end{gather}
\end{example}

\begin{example}\label{ex6}
Note that the underlying directed graphs (unweighted) of $S^nD$ and $S_nD$ 
are the same for $p=\har f=0$, but can be different for $p>0$. 
The diagram of $S_3T$ where $T$ is the standard 2-dimensional representation of 
$sl(2)$ considered in Example \ref{ex5}, is 
\begin{equation}\label{ex6eq1}
\xymatrix{u^3 \ar@(ul,ur)[]^{3h} \ar@/^/[r]^{x_+}
&u^2v \ar@/^/[l]^{-3x_-} \ar@(ul,ur)[]^h \ar@/^/[r]^{2x_+}
&uv^2 \ar@/^/[l]^{-2x_-} \ar@(ul,ur)[]^{-h} \ar@/^/[r]^{3x_+}
&v^3 \ar@(ul,ur)[]^{-3h} \ar@/^/[l]^{-x_-}} .
\end{equation}
The edges labeled by $\pm 3x_\pm$ 
that we have to delete for $p=3$ from the diagrams 
\eqref{ex5eq4} and \eqref{ex6eq1}, have opposite directions.
\end{example}  

\begin{prop}\label{prop10}
If $D$ is a diagram of a representation $T$ of a Lie algebra $L$ in $V$ 
corresponding to the basis $(e_i)_{i\in I}$, then $\Lambda^nD$ is a diagram of 
the subrepresentation $\Lambda_nT$ of $L$ in $\Lambda_nV\subseteq V^{\otimes n}$
where $\Lambda_nV$ is a subspace of $V^{\otimes n}$ with the basis 
$(e_{i_1}\wedge\dots\wedge e_{i_n})_{(i_1<\dots<i_n)\in \Lambda^nI}$ 
where 
\begin{equation}
e_{i_1}\wedge\dots\wedge e_{i_n}=\sum_{s\in S_n}(-1)^{\inv(s)}
e_{s(i_1)}\otimes\dots\otimes e_{s(i_n)}
\end{equation}
where $\inv(s)$ is the number of inversions in $s\in S_n$; corresponding to 
this basis.
\end{prop} 
\begin{proof}Again, it follows directly from the definitions. \end{proof}

\section{Diagram method}\label{sec3}

I describe here the general method of finding normal forms for the orbits 
of strictly triangular representations of Lie algebras, utilizing the diagrams of 
representations. For simplicity, let us suppose that $f$ is a field of 
characteristic 0. 

\begin{defn}\label{defn6}
We'll call the representation $T$ of a Lie algebra $L$ in the space $E$ with 
linearly ordered basis $(e_i)_i\in I$, (strictly) triangular, 
if linear transformations $T(l)$ for all $l\in L$ have (strictly) lower  
triangular matrices in the given basis.
\end{defn}

\begin{prop}\label{prop11}
The representation $T$ is triangular iff for any edge $(i,j)$ of its 
diagram, $i\leq j$. The representation $T$ is strictly triangular iff 
for any edge $(i,j)$ of its diagram, $i<j$.
\end{prop}
\begin{proof}It follows directly from Definition \ref{defn6}.\end{proof}

\begin{defn}\label{defn7}
For a strictly triangular representation $T$, denote $\Gamma(T)$ the group of 
automorphisms of $E$ of the form 
\begin{equation}\label{d7eq1}
\exp T(l)=\sum_{i=0}^{\infty}\frac{1}{i!}T(l)^i
\end{equation}
with $l\in L$.
\end{defn}

\begin{lemma}\label{lem1}
If $D$ is a diagram of a representation $T$ of a Lie algebra $L$ in $E$ 
corresponding to the basis $(e_i)_{i\in I}$, then for any $l\in L$, the 
matrix of $\exp T(tl)$ in the basis $(e_i)_{i\in I}$ is $^tC(D)(l)$, 
the transposed 
walk matrix of the weighted graph $D$, applied to $l$, where 
walk matrix $C(D)$ is an $I\times I$ matrix such that
\begin{equation}\label{l1eq1} 
C(D)_{ij}=\sum_{n=0}^\infty \frac{c_n(D,i\rightarrow j)}{n!}t^n 
\end{equation}
where 
\begin{equation}\label{l1eq2}
c_n(D,i\rightarrow j)=\sum w(i_1,j_1)\dots w(i_n,j_n)
\end{equation}
with summation over all the walks $(i_1,j_1), \dots, (i_n,j_n)$ of length 
$n$ between $i$ and $j$, and $w$ denotes the weight function.
\end{lemma}
\begin{proof}It follows from \eqref{d7eq1} and the definitions of $\exp$ 
and matrix multiplication. \end{proof}

\begin{lemma}\label{lem2}
If for a strictly triangular representation $T$, $l\in L$, $k\in I$ and 
$x\in E$, one has $((\exp T(tl))x)_j=x_j$ for all $j<k$, $t\in f$ and 
$((\exp T(l))x)_k=x_k+c$ with $c\neq 0$, then 
\begin{equation}\label{l2eq1}
((\exp T(tl))x)_k=x_k+ct
\end{equation}
for any $t\in f$. In particular,
\begin{equation}\label{l2eq2}
((\exp T(-\frac{x_k}{c}l))x)_k=0.
\end{equation}
\end{lemma}
\begin{proof}
Since $((\exp T(tl))x)_j-x_i$ for $j<k$ are polynomials of $t$ equal to 0 
for all $t\in f$ and we suppose that $f$ is of characteristic 0, hence infinite, 
all the coefficients of these polynomials must be 0. It means that  
$c_n(D,i\rightarrow j)=0$ for all $j<k$ and $n\geq 1$, where $D$ is the diagram 
of $T$, as usual. Now, 
writing each walk of length $n\geq 2$ from $i$ to $k$ as a composition of a walk of 
length $n-1$ from $i$ to $j$ and the last edge of the original walk, from $j$ to $k$, 
we obtain from \eqref{l1eq2} 
\begin{equation}\label{l2eq3} 
c_n(D, i\rightarrow k)=\sum_{j\in I} c_{n-1}(D,i\rightarrow j) w(j,k)=0 
\end{equation}
for $n\geq 2$. Thus, 
\begin{equation}\label{l2eq4}
((\exp T(tl))x)_k=x_k + \sum_{j\in I}x_j c_1(D,j\rightarrow k)(tl)
=x_k+\sum_{j<k}x_j w(j,k)(l)t.
\end{equation}
Substituting $t=1$, we find 
\begin{equation}\label{l2eq5}
c=\sum_{j<k}x_j w(j,k)(l).
\end{equation}
Equations \eqref{l2eq4} and \eqref{l2eq5} give us \eqref{l2eq1} and 
\eqref{l2eq2} follows from \eqref{l2eq1}.
\end{proof} 

\begin{defn}\label{defn8}
For a linearly ordered basis $(e_i)_{i\in I}$ of $E$, 
introduce the following quasi-ordering on $E$:
\begin{equation}\label{d8eq1}
\sum_{i\in I}x_ie_i \preceq \sum_{i\in I}y_ie_i\quad 
\text{iff}\quad (x_j\neq 0, y_j=0) \Rightarrow \exists i\leq j, (x_i=0, y_i\neq 0).
\end{equation}
For elements 
$x,y\in E$, we'll write $x\prec y$ iff $x\preceq y$ and $y\not\preceq x$.
\end{defn}

In other words, $x\prec y$ means that in the first place where $x_iy_i=0$ and 
coefficients $x_i, y_i$ are not both $0$, one finds $x_i=0, y_i\neq 0$. If 
$x_iy_i=0$ implies $x_i=y_i=0$, then $x\sim y$, i.\ e.\ $x\preceq y$ and 
$y\preceq x$. 

\begin{theorem}\label{thm1}
For a strictly triangular representation $T$ in a finite dimensional $f$-space $E$, 
every orbit of $\Gamma(T)$ in $E$ 
has the unique lowest point with respect to $\prec$; i.\ e., the point 
$x$ such that $x\prec y$ for any $y\in \Gamma(T)x, y\neq x$. 
\end{theorem}
\begin{proof}
Choose an arbitrary $\Gamma(T)$-orbit and a point $x^{(0)}$ on it. Look at the 
coordinates of this point according to the order of $I$, trying to turn non-zero 
coordinates into 0 by applying elements of $\Gamma(T)$ to $x^{(0)}$. Since $T$ is
strictly triangular, the action of $\Gamma(T)$ on each coordinate $x^{(0)}_j$ 
depends 
only on the preceding coordinates $x^{(0)}_i$, with $i<j$. In particular, we can't 
change the first non-zero coordinate. Let $i_0$ be the smallest index of a 
non-zero coordinate of $x^{(0)}$ that we can change, leaving all the preceding 0's as 
0's. Then, applying Lemma \ref{lem2}, we can make $(\alpha_1x^{(0)})_{i_0}=0$  
by action of some element $\alpha_1\in\Gamma(T)$, leaving all the preceding 
0's as 0's. Denote $x^{(1)}=
\alpha_1x^{(0)}$. The element $x^{(1)}$ has the same coordinates 
$x^{(1)}_i=x^{(0)}_i$ as the element $x^{(0)}$ for $i<i_0$, by construction, 
and $x^{(1)}_{i_0}=0$ while $x^{(0)}_{i_0}\neq0$. Thus, $x^{(1)}\prec x^{(0)}$. 
Analagously, denoting $i_1$ the smallest index of a non-zero coordinate of 
$x^{(1)}$ that we can change, leaving all the preceding 0's as 0's, applying 
Lemma \ref{lem2}, we can 
make $(\alpha_2x^{(1)})_{i_1}=0$  
by action of some element $\alpha_1\in\Gamma(T)$, leaving all the preceding 
0's as 0's. Denote $x^{(2)}=\alpha_2x^{(1)}$. The same as above, we have 
$x^{(2)}\prec x^{(1)}\prec x^{(0)}$. Continuing this construction for a sequence of indices 
$i_0<i_1<\dots<i_k$, we can construct a sequence of points 
$x^{(k+1)}\prec\dots\prec x^{(0)}$ of the $\Gamma(T)$-orbit of $x^{(0)}$. 
Because of the finiteness of $I$, earlier or later, we must stop on the point 
of $x^{(k+1)}$ of this sequence. By our construction, we can't change any non-zero 
coordinates of $x^{(k+1)}$ leaving all the preceding 0's as 0's, by elements 
of $\Gamma(T)$. It means that $x^{(k+1)}$ is $\prec$ than any other point on the 
orbit. We constructed the lowest point on the orbit. 
\end{proof} 

\begin{defn}\label{defn9}
We'll call the lowest point on an orbit with respect to $\prec$, 
by the normal form of the orbit. 
\end{defn}

\begin{lemma}\label{lem3}
If for a strictly triangular representation $T$, $l\in L$, $k\in I$ and 
$x\in E$, one has $(T(l)x)_j=0$ for all $j<k$ and 
$(T(l)x)_k=c\neq 0$, then $((\exp T(tl))x)_j=x_j$ for all $j<k$, $t\in f$ and
\begin{equation}\label{l3eq1}
((\exp T(tl))x)_k=x_k+ct
\end{equation}
for any $t\in f$. In particular,
\begin{equation}\label{l3eq2}
((\exp T(-\frac{x_k}{c}l))x)_k=0.
\end{equation}
\end{lemma}
\begin{proof}Note that $(T(l)x)_j=0$ implies $(T(tl)x)_j=0$. Now, because $T$ is 
strictly triangular and $(T(tl)x)_j=0$ for all $j<k$, we have 
$(T(tl)^2x)_j=0$ for all $j\leq k$. Thus 
\begin{equation}\label{l3eq4}
((\exp T(tl))x)_j=((1+T(tl))x)_j=x_j
\end{equation}
for $j<k$ and
\begin{equation}\label{l3eq3}
((\exp T(tl))x)_k=((1+T(tl))x)_k=x_k+ct ,
\end{equation}
q.e.d. Formula \eqref{l3eq2} follows from \eqref{l3eq1}.\end{proof}

\begin{theorem}\label{thm2}
For a strictly triangular representation $T$ of a Lie algebra $L$ 
in a finite dimensional $f$-space $E$,
the element $x\in E$ is a normal form on its $(\exp T)$-orbit, iff 
\begin{equation}\label{t2eq1}
x\preceq x+T(l)x  
\end{equation}
for all $l\in L$.
\end{theorem}
\begin{proof}
Proof by contradiction. Suppose that an $x\in E$ is a normal form on its orbit, 
and the set of $l\in L$ such that
\begin{equation}\label{t2eq2}
x+T(l)x\prec x ,
\end{equation}
is not empty, i.\ e.\ there exist such $k\in I$, $l\in L$ that 
$x_k\neq 0$ but $(x+T(l)x)=0$ and $(T(l)x)_j=0$ if $j<k$ and $x_j=0$. 
In the last case $x_k\neq 0$, but
\begin{equation}\label{t2eq3}
(T(l)x)_k=(x+T(l)x)_k-x_k=-x_k\neq 0.
\end{equation}
Let $k\in I$ be the smallest, for all $l\in L$, index such that $x_k\neq 0$ and 
$(T(l)x)_k\neq 0$ and $(T(l)x)_j=0$ if $j<k$ and $x_j=0$. Then 
$(T(l)x)_j=0$ for all $j<k$. Applying Lemma \ref{lem3}, we get 
\begin{equation}\label{t2eq4}
((\exp T(-\frac{x_k}{c}l))x)_j=\begin{cases}x_j &\text{for $j<k$,}\\
0 &\text{for $j=k$,}\end{cases}
\end{equation}
that means
\begin{equation}\label{t2eq5}
(\exp T(-\frac{x_k}{c}l))x\prec x\thinspace ;
\end{equation}
a contradiction. Thus, for every normal form we have \eqref{t2eq1} for all $l\in L$.
Conversely, if for some $x\in E$ we have \eqref{t2eq1} for all $l\in L$, then 
let $k\in I$ be the smallest, for all $l\in L$, index such that $x_k\neq 0$ and 
$(\exp T(l)x)_k=x_k+c$ with $c\neq 0$ and $(\exp T(l)x)_j=0$ if $j<k$ and $x_j=0$. 
Applying Lemma \ref{lem2}, we get 
\begin{equation}\label{t2eq6}
x+T(-\frac{x_k}{c}l)x\prec x, 
\end{equation}
that is a contradiction to \eqref{t2eq1}. Thus 
if for some $x\in E$ we have \eqref{t2eq1} for all $l\in L$, then $x$ is a normal 
form on its $\exp(T)$-orbit.
\end{proof}

\begin{cor}\label{cor2}
For a strictly triangular representation $T$ of a Lie algebra $L$ 
in a finite dimensional $f$-space $E$,
the element $x\in E$ is a normal form on its $(\exp T)$-orbit, iff 
for every $k\in I$ with $x_k\neq 0$, 
\begin{equation}\label{c2eq1}
(T(l)x)_k\subseteq \spa\{((T(l)x)_i)_{i<k}\}\subseteq L^\ast . 
\end{equation}
\end{cor}
\begin{proof}
The condition \eqref{c2eq1} means that we can't change non-zero coordinates 
of $x$ by adding $T(l)x$ without changing the previous coordinates, which is 
the same as the stament of Theorem \ref{thm2}.
\end{proof}

The diagram $D$ of a representation $L$ gives us a convenient way of writing 
down the elements 
\begin{equation}\label{c2eq2}
(T(l)x)_j=\sum_{i<j}x_iw(i,j)(l)
\end{equation}
for a strictly triangular representation $T$, see \eqref{eqpr6} and Proposition 
\ref{prop11}. 

Using diagrams and Corollary \ref{cor2}, 
we can construct all the normal forms of strictly triangular 
representations. Let $T$ be a strictly triangular representation of a Lie algebra 
$L$ in the $f$-space $E$ with linearly ordered basis $(e_i)_{1\leq i\leq n}$. 
Denote $D$ the diagram of $T$ and $w$ its weight function. 
First we construct 
the normal form `in a general position'. Put $x_1=c_1\in f^\ast$ with an arbitrary 
$c_1\in f^\ast$. If $w(1,2)\neq 0$, put $x_2=0$ and remember the element 
$y_2=c_1w(1,2)\in L^\ast$, if $w(1,3)\not\in y_2f\subseteq L^\ast$, put $x_3=0$ 
and remember $y_3=c_1w(1,3)\in L^\ast$, \ldots, until $w(1,k)\in \spa\{y_2, \dots, 
y_{k-1}\}\subseteq L^\ast$, in which case put $x_k=c_k\in f^\ast$; now if 
$c_1w_{1,k+1}+c_kw_{k,k+1}\not\in \spa\{y_2, \dots, y_{k-1}\}\subseteq L^\ast$, put 
$x_{k+1}=0$ and remember $y_{k+1}=c_1w_{1,k+1}+c_kw_{k,k+1} \in L^\ast$, otherwise 
put $x_{k+1}=c_{k+1}\in f^\ast$ and don't remember the value of $y_{k+1}$ and so on, 
assigning $x_j=0$ and memorizing the corresponding value of \eqref{c2eq2}, denoting 
it $y_j$ if it does not belong to the $\spa\{y_i\}_{i<j}$ and assigning $x_j=c_j
\in f^\ast$ otherwise, without remembering $y_j$ in that case, until we define 
the value of the last coefficient, $x_n$.

\begin{theorem}\label{thm3}
The element 
\begin{equation}\label{t3eq1}
x=\sum_{i\in I} x_ie_i \thickspace\in E
\end{equation}
with coefficients $x_i$ described above, is a normal form.
\end{theorem}
\begin{proof}
It follows from Corollary \ref{cor2}, formula \eqref{c2eq2} and our construction. 
\end{proof}

\begin{defn}\label{defn10}
Let us call the element $x\in E$ of form \eqref{t3eq1} with $x_i$ given by the 
construction above, a normal form in a general position.
\end{defn}

\begin{example}\label{ex7}
The adjoint representation of the Lie algebra of strictly upper triangular 
$4\times 4$ matrices, has a diagram
\begin{equation}\label{ex7eq1}
\xymatrix{x_4 \ar[r]_{-a_{23}} \ar@/^1pc/[rr]^{-a_{24}}
&x_5 \ar[r]_{-a_{34}}&x_6\\
&x_2 \ar[r]_{-a_{34}} \ar[u]^{a_{12}}&x_3 \ar[u]^{a_{12}}\\
&&x_1 \ar[u]^{a_{23}} \ar@/_1pc/[uu]_{a_{13}}}.
\end{equation}
Applying the diagram method, we get a normal form in a general position
\begin{equation}\label{ex7eq2}
\begin{pmatrix}&c_4&0&0\\&&c_2&0\\&&&c_1\\ &&& \end{pmatrix}
\end{equation}
with non-zero $c_1, c_2, c_4$. The memorized linear forms are 
$y_3=c_1a_{23}-c_2a_{34}, \thickspace y_5=c_2a_{12}-c_4a_{23}$ and 
$y_6=c_1a_{13}-c_4a_{24}$. Also, by Proposition \ref{prop1}, the co-adjoint 
representation of this Lie algebra has a diagram
\begin{equation}\label{ex7eq3}
\xymatrix{x_3\\ x_2 \ar[u]_{a_{23}} \ar[r]^{-a_{12}} &x_5\\
x_1 \ar@/^1pc/[uu]^{a_{24}} \ar[u]_{a_{34}} \ar[r]^{-a_{12}} \ar@/_1pc/[rr]_{-a_{13}}
&x_4 \ar[u]_{a_{34}} \ar[r]^{-a_{23}}&x_6}.
\end{equation}
Applying the diagram method, we get a normal form in a general position
\begin{equation}\label{ex7eq4}
\begin{pmatrix}&&&\\ 0&&&\\ 0&c_5&&\\ c_1&0&0& \end{pmatrix} 
\end{equation}
with non-zero $c_1, c_5$.
\end{example}

After describing the normal forms in a general position, we can continue 
the description of normal forms as follows. Take a normal form $x$, put 
$x_k=0$ instead of the last non-zero coefficient of it and continue the procedure 
as above, for $x_{k+1}$ and so on, the same as we did calculating 
the normal forms in a general position, saving all the memorized linear forms
$y_i$ for $i<k$ and changing them appropriately for $i>k$, 
until we define the value of the last coefficient, $x_n$.

\begin{theorem}\label{thm4}
If we apply the procedure described in the previous paragraph to a normal form, 
we'll get a normal form as well.
\end{theorem}
\begin{proof}
Indeed, all the conditions of Corollary \ref{cor2} are true for $i<k$ because 
we started from a normal form and they are also true for $i\geq k$, by construction. 
Thus, by Corollary \ref{cor2}, we get a normal form.\end{proof}

Then, applying the same procedure to a new normal form, we can construct another 
one and so on until we get $0$, where we stop.
This method, I called a diagram method because of the intensive use of 
diagrams in it, is not an algorithm for an infinite field in a rigorous sence 
of this word because it sometimes requires an infinite number of steps. 
Nevertheless, it allows us to describe the normal forms in many cases rather easily. 
A few years ago I wrote down all the normal 
forms for the co-adjoint representations of the Lie algebras of strictly upper 
triangular $n\times n$ matrices with $n\leq 8$ without using a computer. 
For a finite field, it is a real algorithm. We can't use the exponential 
formula in that case in general, but we can use another correspondence between 
algebraic Lie groups and algebras. For the case of the Lie algebra of strictly 
triangular matrices over a finite field, one can use the correspondence 
$T\rightarrow 1+T$, see the details in \cite{KF}. 

\begin{example}\label{ex8}
Applying the diagram method to the normal formal in a general position, 
\eqref{ex7eq4}, we obtain the following complete list of normal forms 
for that case:
\begin{equation}\label{ex8eq1}\begin{split}
&\begin{pmatrix}&&&\\ 0&&&\\ 0&c_5&&\\ c_1&0&0& \end{pmatrix}\rightarrow 
\begin{pmatrix}&&&\\ 0&&&\\ 0&0&&\\ c_1&0&0& \end{pmatrix}\rightarrow 
\begin{pmatrix}&&&\\ 0&&&\\ c_2&0&&\\ 0&c_4&c_6& \end{pmatrix}\rightarrow
\begin{pmatrix}&&&\\ 0&&&\\ c_2&0&&\\ 0&c_4&0& \end{pmatrix}\rightarrow\\ 
&\begin{pmatrix}&&&\\ 0&&&\\ c_2&0&&\\ 0&0&c_6& \end{pmatrix}\rightarrow
\begin{pmatrix}&&&\\ 0&&&\\ c_2&0&&\\ 0&0&0& \end{pmatrix}\rightarrow
\begin{pmatrix}&&&\\ c_3&&&\\ 0&0&&\\ 0&c_4&0& \end{pmatrix}\rightarrow
\begin{pmatrix}&&&\\ c_3&&&\\ 0&c_5&&\\ 0&0&c_6& \end{pmatrix}\rightarrow\\  
&\begin{pmatrix}&&&\\ c_3&&&\\ 0&c_5&&\\ 0&0&0& \end{pmatrix}\rightarrow
\begin{pmatrix}&&&\\ c_3&&&\\ 0&0&&\\ 0&0&c_6& \end{pmatrix}\rightarrow
\begin{pmatrix}&&&\\ c_3&&&\\ 0&0&&\\ 0&0&0& \end{pmatrix}\rightarrow
\begin{pmatrix}&&&\\ 0&&&\\ 0&0&&\\ 0&c_4&0& \end{pmatrix}\rightarrow\\
&\begin{pmatrix}&&&\\ 0&&&\\ 0&c_5&&\\ 0&0&c_6& \end{pmatrix}\rightarrow
\begin{pmatrix}&&&\\ 0&&&\\ 0&c_5&&\\ 0&0&0& \end{pmatrix}\rightarrow
\begin{pmatrix}&&&\\ 0&&&\\ 0&0&&\\ 0&0&c_6& \end{pmatrix}\rightarrow
\begin{pmatrix}&&&\\ 0&&&\\ 0&0&&\\ 0&0&0& \end{pmatrix}\end{split}
\end{equation}
with non-zero coefficients $c_i$.
\end{example}

Notice that for the second matrix in this example, the set of 
memorized linear forms $(y_k)$ is exactly the same as for the first matrix which 
is a normal form 
in a general position. In such cases, when the first normal forms 
obtained from a normal form in a general position have the same set of the 
memorized $y_k$, it is convenient to refer to these new normal forms 
as forms in a general position as well; and we'll do that in the following two 
sections. 

Normal forms in the last example form a few families parametrized as $(f^\ast)^k$, 
without any additional conditions. The same is true for the co-adjoint 
representations of Lie algebras of strictly upper triangular $n\times n$ 
matrices with $n\leq 8$. 
However, the following example shows that for $n\geq 9$ there are 
some additional polynomial conditions. 

\begin{prop}\label{prop12}
The matrix
\begin{equation}\label{p12eq1}
\begin{pmatrix}&&&&&&&&\\
0&&&&&&&&\\
0&0&&&&&&&\\
0&0&0&&&&&&\\
0&c_{13}&0&0&&&&&\\
c_4&0&0&0&0&&&&\\
0&0&c_{18}&0&0&x_{33}&&&\\
0&0&0&0&c_{28}&c_{32}&x_{35}&&\\
0&0&0&c_{22}&c_{27}&c_{31}&0&0&
\end{pmatrix}
\end{equation}
with $c_i\in f^\ast$ and at least one of the coefficients $x_{33}$ and 
$x_{35}$ being non-zero, is a normal form for a co-adjoint representation of a Lie algebra of 
strictly upper-triangular $9\times 9$ matrices, iff 
\begin{equation}\label{p12eq2}
\begin{vmatrix}c_{28}&c_{32}\\ c_{27}&c_{31} \end{vmatrix} = 
c_{28}c_{31}-c_{27}c_{32} =0 .
\end{equation}
\end{prop}
\begin{proof}Only if the determinant \eqref{p12eq2} is 0, the corresponding 
linear forms for $x_{33}$ or $x_{35}$ are linearly dependent on the memorized 
earlier linear forms, corresponding to the diagram method.
\end{proof}

\section{Normal forms of quadratic differentials on a line}\label{sec4}

Let $f$ be a field of characteristic $0$. Denote $L_n$ the Lie algebra generated 
by $(l_i)_{i\geq n}$ 
with $i\in \ZZ$, $n>0$ and 
\begin{equation}\label{s2eq1}
[l_i,l_j]=(j-i)l_{i+j}, 
\end{equation}
see \cite{MF}. Denote $G_n$ the group of automorphisms $g$ of the algebra 
$f[[t]]$ of formal power series, such that $g(t)=t+o(t^n)$. 
$L_n$ can be considered as a Lie algebra of $G_n$, with 
\begin{equation}\label{s2eq2}
l_i=t^{i+1}\frac{d}{dt}. 
\end{equation}
Denote also 
$L_{mn}=L_m/L_n$ and $G_{mn}=G_m/G_n$ for $m\leq n$. 

As it was explained in \cite{MF}, we can realize $L_n^\ast$ as a space of 
quadratic differentials with the basis $(y_i=t^{-i-2}(dt)^2)_{i\geq n}$ dual to 
$(l_i)_{i\geq n}$ and $L_{mn}^\ast$ as a subspace of $L_m^\ast$ with the basis 
$(y_i)_{m\leq i<n}$. 

\begin{theorem}\label{s2thm1}
Let $m\geq 0$. If $n\leq 2m+2$, then all the elements of $L_{m+1,n+1}^\ast$ 
are normal 
forms in a general position. If $n>2m+2$, normal forms in a general position for 
the coadjoint representation of $L_{m+1,n+1}$, are $c_0y_n+c_1y_{n-1}+ \dots +
c_my_{n-m}$ with any $c_0\in f^\ast, c_1, \dots, c_m \in f$, for odd $n$, or 
$c_0y_n+c_1y_{n-1}+ \dots +c_my_{n-m}+c_{n/2}y_{n/2}$ 
with any $c_0\in f^\ast, c_1, \dots, c_m, c_{n/2} \in f$, for even $n$. Each normal 
form in a general position in $L^\ast_{m+1,n+1}$ is a normal form in 
$L^\ast_{m+1,i+1}$ for each $i\geq n$ and in $L^\ast_{m+1}$. Each normal form 
in $L^\ast_{m+1,i+1}$ is a normal form in a general position in an  
$L^\ast_{m+1,n+1}$ with $\min\{2m+2,i\}\leq n \leq i$. Each normal form 
in $L^\ast_{m+1}$ is a normal form in a general position in an 
$L^\ast_{m+1,n+1}$ with $n\geq 2m+2$.
\end{theorem}
\begin{proof}
The diagrams of the co-adjoint representations of $L_{m+1,n+1}$ and of $L_{m+1}$
have an edge $(k,j)$ iff $k>j+m$ and $k\neq 2j$ and 
\begin{equation}\label{s2t1eq1}
w(k,j)=(k-2j)y_{k-j}
\end{equation}
in this case. Applying the diagram method, we immediately get the result of 
Theorem \ref{s2thm1}. 
\end{proof}

\begin{cor}\label{s2cor1}
For odd $n$, the ring of polynomial invariants of the co-adjoint representation 
of $L_{m+1,n+1}$ is $f[l_n,\dots,l_{n-m}]$. For even $n$, this ring is 
$f[l_n,\dots,l_{n-m}, P]$ where $P$ is a polynomial.
\end{cor}
\begin{proof}
Rational invariants of a representation give us the equations 
parametrizing the orbits in a general position. Computing the orbits 
of normal forms in a general position, given by Theorem \ref{s2thm1}, we can see 
that they are affine subspaces $l_n=c_0, \dots, l_{n-m}=c_m$ for odd $n$, or, 
for even $n$,  
intersections of the affine subspaces given by the same equations and an 
affine hypersurface, the equation of which must be given by a polynomial, by 
Luroth theorem.\end{proof} 

Theorem \ref{s2thm1} and Corollary \ref{s2cor1} belong to Alexandre Kirillov 
\cite{KV}. The explicit formulas for the polynomials $P$ were found in my 
work \cite{MF}.

\section{Normal forms of generalized formal tensor fields on a line}\label{sec5}

Let $f$ be a field of characteristis $0$, Lie algebras $L_{m+1}$ and groups 
$G_{m+1}$ defined in the previous section. Denote $F_{\lambda\mu} = 
f[[t]]t^\mu(dt)^{-\lambda}$. Lie algebras $L_{m+1}$ and groups $G_{m+1}$ 
naturally act on $F_{\lambda\mu}$. We'll use the topological basis 
$(e_n=t^{n+\mu}(dt)^{-\lambda})_{n\in \ZZ_{\geq 0}}$ of $F_{\lambda\mu}$.

\begin{theorem}\label{s3thm1}
Let $m\geq 0$. If $\mu\neq (m+k+1)\lambda$ with a positive integer $k$, then the 
normal forms in a general position in $F_{\lambda\mu}$, are $c_0e_0 + 
c_1e_1 + \dots +c_me_m$ with any $c_0\in f^\ast, c_1, \dots, c_m\in f$. 
If $\lambda\neq 0$ and $\mu=(m+k+1)\lambda$ for a positive integer $k$, the 
normal forms in general position in $F_{\lambda\mu}$, are 
$c_0e_0 + c_1e_1 + \dots +c_me_m + c_{m+k}e_{m+k}$ 
with any $c_0\in f^\ast, c_1, \dots, c_m, c_{m+k}\in f$. Each normal form in a 
general position in $F_{\lambda\mu}$, is a normal form in $F_{\lambda,\mu-n}$ 
for every nonnegative integer $n$.  Each normal form 
in $F_{\lambda\mu}$, is a normal form in a general position in $F_{\lambda,\mu+n}$ 
for a nonnegative integer $n$; except for the case $\lambda=0$ and $\mu$ is 
a non-positive integer, in which case there are additional normal forms, 
$c+C$, for any constant $c\in f^\ast$ and $C$, a normal form in a general position 
in $F_{0n}$ with a positive integer $n$.
\end{theorem}
\begin{proof}The diagram of the representation of $L_{m+1}$ in     
$F_{\lambda\mu}$ has an edge $(j,k)$ iff $k>j+m$ and $j+\mu\neq (k-j+1)\lambda$ ,
in which case 
\begin{equation}\label{s3t1eq1}
w(j,k)=(j+\mu-(k-j+1)\lambda)y_{k-j}.
\end{equation}
Applying the diagram method, we immediately get the result of Theorem 
\ref{s3thm1}.\end{proof} 

Analogously Corollary \ref{s2cor1}, we can describe the invariants of the 
representation of $L_{m+1}$ in $F_\lambda\mu$, see the explicit formulas in 
\cite{MF}.

\begin{acn}
I would like to thank Alexandre Kirillov and Fan Chung Graham, as well as my mother  
and my Beautiful and Wonderful wife, Bette.
\end{acn}

\end{document}